\documentclass[a4paper,11pt]{article}
\title{\bf{On a certain generalization of spherical twists}
}
\date{}
\author{Yukinobu Toda}
\usepackage{amsmath,amssymb}
\usepackage{macros}

\usepackage{array}
\usepackage{amscd}
\usepackage[all]{xy}

\DeclareFontFamily{U}{rsfs}{%
\skewchar\font127}
\DeclareFontShape{U}{rsfs}{m}{n}{%
<-6>rsfs5<6-8.5>rsfs7<8.5->rsfs10}{}
\DeclareSymbolFont{rsfs}{U}{rsfs}{m}{n}
\DeclareSymbolFontAlphabet
{\mathrsfs}{rsfs}
\DeclareRobustCommand*\rsfs{%
\@fontswitch\relax\mathrsfs}

\newtheorem{thm}{Theorem}[section]
\newtheorem{prop}[thm]{Proposition}

\newtheorem{defi}[thm]{Definition}
\newtheorem{rmk}[thm]{Remark}

\newtheorem{prop-defi}[thm]{Proposition-Definition}

\begin{document}
\maketitle
\begin{abstract}
This note gives a generalization of spherical twists, 
and describe the autoequivalences associated to certain
non-spherical objects.   
Typically these
 are obtained by deforming the 
 structure sheaves of
$(0, -2)$-curves on threefolds, or deforming 
$\mathbb{P}$-objects introduced by 
D.Huybrecht and R.Thomas. 
\end{abstract} 
\section{Introduction}
In this paper, we introduce a new class of autoequivalences 
of derived categories of coherent sheaves
on smooth projective varieties, which generalizes the notion of
spherical twists given in~\cite{ST}. Such autoequivalences 
are associated to a certain class of objects, which are not necessary
spherical but 
are interpreted as ``fat" version of them. 
We introduce the notion of $R$-spherical objects for a noetherian and artinian 
local $\mathbb{C}$-algebra $R$, and imitate the construction of 
spherical twists to give the associated autoequivalences. 

\vspace{3mm}

Let $X$ be a smooth complex projective variety, 
and $D(X)$ be a bounded derived category of coherent sheaves 
on $X$. When $X$ is a Calabi-Yau 
3-fold, $D(X)$ is considered to 
represent the category of $D$-branes of type $B$, and should
be equivalent to the derived Fukaya category on a mirror 
manifold under Homological mirror symmetry~\cite{Kon}.   
On the mirror side, there are typical symplectic automorphisms
by taking Dehn twists along Lagrangian spheres. 
The 
notions of spherical objects and associated twists were 
introduced in~\cite{ST} in order
 to realize Dehn twists
under mirror symmetry.  
Recall that $E\in D(X)$ is called spherical if and only if the 
following holds~\cite{ST}:
\begin{itemize}
\item $\Ext _X ^i (E, E)= \left\{ \begin{array}{cl}
 \mathbb{C} & (i=0, \dim X), \\ 0 & (\mbox{otherwise}).   \end{array} \right.  $\item $E \otimes \omega _X \cong E .$
 \end{itemize}
Then one can construct the autoequivalence $T_{E}\colon D(X)\to D(X)$ 
which fits into the distinguished triangle~\cite{ST}:
$$\dR \Hom (E, F)\otimes _{\mathbb{C}} E \lr F \lr T_E (F), $$
for $F\in D(X)$.  $T_E$ is called a spherical twist. 
This is a particularly important class of autoequivalences, 
especially when we consider $A_n$-configulations on surfaces 
 as indicated in~\cite{UI}. 
On the other hand, it has been observed that 
there are some autoequivalences which are not described in terms of 
spherical twists.
This occurs even in the similar situation discussed in~\cite{UI} as 
follows.   
Let $X\to Y$ be a three dimensional 
flopping contraction which contracts a rational 
curve $C\subset X$, and $X^{\dag}\to Y$ be its flop. 
Then one can construct the autoequivalence~\cite{B-O2},~\cite{Br1},~\cite{Ch}, 
$$\Phi \cneq 
\Phi _{X^{\dag}\to X}^{\oO _{X \times _Y X^{\dag}}}
\circ \Phi _{X \to X^{\dag}}^{\oO _{X \times _Y X^{\dag}}}
\colon D(X) \lr D(X^{\dag}) \lr D(X). $$
If $C\subset X$ is not a $(-1, -1)$-curve, $\Phi$ is not written as 
a spherical twist, and our motivation comes from describing such 
autoequivalences. 
Let $R$ be a noetherian and 
 artinian local $\mathbb{C}$-algebra. 
We introduce the notion of $R$-spherical objects defined on $D(X\times 
\Spec R)$. In the above example, $\Spec R$ is taken to be the moduli space of 
$\oO _C (-1)$, and the universal family gives
the $R$-spherical object.  Our main theorem is the following:
\begin{thm}
To any $R$-spherical object $\eE \in D(X\times \Spec R)$, we can 
associate the autoequivalence $T_{\eE}\colon D(X)\to D(X)$, which 
fits into the distinguished triangle:
$$\dR \Hom _X(\pi _{\ast}\eE, F) \dotimes _R \pi _{\ast}\eE \lr 
F \lr T_{\eE}(F), $$ 
for $F\in D(X)$. Here $\pi \colon X\times \Spec R\to X$ is the projection. 
\end{thm}
Using the notion of $R$-spherical objects and associated twists, 
we can also give the deformations of $\mathbb{P}$-twists
in the case which is not treated in~\cite{HT}. 

\subsection*{Acknowledgement}
The author thanks Tom Bridgeland for useful discussions and comments. 
He is supported by Japan Society for the Promotion of Sciences Research 
Fellowships for Young Scientists, No 1611452.

\subsection*{Notations and conventions}
\begin{itemize}
\item For a variety $X$, we denote by $D(X)$ its bounded derived category 
of coherent sheaves. 
\item $\Delta$ means the diagonal $\Delta \subset X\times X$ or 
the diagonal embedding $\Delta \colon X\to X\times X$. 
\item For another variety $Y$ and an object $\pP \in D(X\times Y)$, 
denote by $\Phi _{X\to Y}^{\pP}$ the integral transform with kernel 
$\pP$, i.e. 
$$\Phi _{X\to Y}^{\pP}(\ast)\cneq \dR p_{Y\ast}(p_X ^{\ast}(\ast)\dotimes 
\pP) \colon D(X) \lr D(Y).$$
Here $p_X$, $p_Y$ are projections from $X\times Y$ onto corresponding 
factors.
\end{itemize}

\section{Generalized spherical twists}
Let $X$ be a smooth projective variety over $\mathbb{C}$ and
$R$ be a noetherian and 
 artinian local $\mathbb{C}$-algebra. 
We introduce the notion of $R$-spherical objects defined on  
$D(X\times \Spec R)$. 
Let $\pi \colon X\times \Spec R \to X$ and $\pi ' \colon X\times \Spec R 
\to \Spec R$ be 
projections and $0\in \Spec R$ be the closed point. 
\begin{defi}
Take $\eE \in D(X\times \Spec R)$ and let 
$E\cneq \eE |_{X\times \{ 0\}} \in D(X)$
be its derived restriction. Then $\eE$ is $R$-spherical if and only if the 
following conditions hold:
\begin{itemize}
\item $\Ext _X ^i (\pi _{\ast}\eE, E)= \left\{ \begin{array}{cl}
 \mathbb{C} & (i=0, \dim X), \\ 0 & (\emph{otherwise}).   \end{array} \right.  $\item $E \otimes \omega _X \cong E .$
 \end{itemize}
\end{defi}
\begin{rmk}
If $R=\mathbb{C}$, then $R$-spherical objects coincide with 
 usual spherical objects.
\end{rmk}
We imitate the construction of the spherical twists in the 
following theorem. 
\begin{thm}\label{equ}
To any $R$-spherical object $\eE \in D(X\times \Spec R)$, we can 
associate the autoequivalence $T_{\eE}\colon D(X)\to D(X)$, which 
fits into the distinguished triangle:
$$\dR \Hom _X(\pi _{\ast}\eE, F) \dotimes _R \pi _{\ast}\eE \lr 
F \lr T_{\eE}(F), $$ 
for $F\in D(X)$. 
\end{thm}
\textit{Proof}. 
First we construct the kernel of $T_{\eE}$. Let 
$p_{ij}$ and $p_i$ 
be projections as in the following diagram:
$$
\xymatrix{
& & X\times \Spec R \times X 
\ar[dl]_{p_{12}}
\ar[d]^{p_{13}}
\ar[dr]^{p_{23}}
& & \\
 & X\times \Spec R \ar[d]^{\pi} &  X\times X \ar[dl]_{p_1}\ar[dr]^{p_2}
  & X\times \Spec R \ar[d]^{\pi} , \\
 & X & & X &
}
 $$
and consider the object,
$$\qQ \cneq 
\dR p_{13\ast}\left( p_{12}^{\ast}(\pi ^{!}\oO _X \dotimes \check{\eE} )
\dotimes p_{23}^{\ast}\eE \right) \in D(X\times X). $$
Here $\check{\eE}
$ means its derived dual. 
Then for $F\in D(X)$, we can calculate $\Phi _{X\to X}^{\qQ}(F)$ as follows:
\begin{align*}
\Phi _{X\to X}^{\qQ}(F) &\cong \dR p_{2\ast}(\dR p_{13\ast}(p_{12}^{\ast}
(\pi ^{!}\oO _X \dotimes \check{\eE})\dotimes p_{23}^{\ast}\eE) \dotimes 
p_1 ^{\ast}F) \\
&\cong \dR p_{2\ast}\dR p_{13\ast}(p_{12}^{\ast}(\pi ^{!}\oO _X \dotimes
\check{\eE})\dotimes p_{23}^{\ast}\eE \dotimes p_{13}^{\ast}p_1 ^{\ast}F)
\\ 
&\cong \pi _{\ast}\dR p_{23\ast}(p_{12}^{\ast}(\pi ^{!}\oO _X \dotimes 
\check{\eE})\dotimes p_{23}^{\ast}\eE \dotimes p_{12}^{\ast}\pi ^{\ast}F) \\
&\cong 
\pi _{\ast}\left\{ \eE \dotimes \dR p_{23\ast}p_{12}^{\ast}(\pi ^{!}\oO _X 
\dotimes \check{\eE} \dotimes \pi ^{\ast}F) \right\} \\
&\cong 
\pi _{\ast} \left\{ \eE \dotimes \pi ^{'\ast}\dR \pi ' _{\ast}\dR 
\hH om (\eE, \pi ^{!}F)\right\} \\
&\cong \pi _{\ast}\eE \dotimes _R \dR \Hom (\pi _{\ast}\eE, F)  . 
\end{align*}
The fifth equality comes from the base change formula for the diagram below:
$$\xymatrix{
X\times \Spec R \times X \ar[r]^{p_{12}}\ar[d]_{p_{23}}
 & X\times \Spec R \ar[d]^{\pi '} \\
X\times \Spec R \ar[r]_{\pi '} & \Spec R .
}$$ 
On the other hand, we have 
\begin{align*}
\Hom _{X\times X}(\qQ, \oO _{\Delta}) &= \Hom _{X\times X}
(\dR p_{13\ast}( p_{12}^{\ast}(\pi ^{!}\oO _X \dotimes \check{\eE} )
\dotimes p_{23}^{\ast}\eE ), \oO _{\Delta}) \\
&= \Hom _X(\dL \Delta ^{\ast}\dR p_{13\ast}( p_{12}^{\ast}(\pi ^{!}\oO _X \dotimes \check{\eE} )
\dotimes p_{23}^{\ast}\eE ), \oO _{X}) \\
&= \Hom _X (\pi _{\ast}\dL (\Delta, \id)^{\ast}( p_{12}^{\ast}(\pi ^{!}\oO _X \dotimes \check{\eE} )
\dotimes p_{23}^{\ast}\eE ), \oO _{X}) \\
&= \Hom _X (\dL (\Delta, \id)^{\ast}( p_{12}^{\ast}(\pi ^{!}\oO _X \dotimes \check{\eE} )
\dotimes p_{23}^{\ast}\eE ), \pi ^{!}\oO _{X}) \\
&= \Hom _X(\pi ^{!}\oO _X \dotimes \check{\eE}\dotimes \eE, \pi ^{!}\oO _X). 
\end{align*}
The third equality comes from the base change formula for the diagram below:
$$
\xymatrix{
X \times \Spec R \ar[r]^{(\Delta, \id)} 
\ar[d]_{\pi} & X\times \Spec R \times X \ar[d]^{p_{13}} \\
X \ar[r]_{\Delta} & X\times X . \\
}$$
Let $\mu \colon \qQ \to \oO _{\Delta}$ be the morphism which
 corresponds to the morphism
$$\id _{\pi ^{!}\oO _X}\otimes \ev \colon
\pi ^{!}\oO _X \dotimes \check{\eE}\dotimes \eE \lr \pi ^{!}\oO _X ,$$
under the above isomorphisms. Let us take its cone 
$\rR \cneq 
\CCone (\mu)\in D(X \times X)$. Then the above 
calculation for $\Phi _{X\to X}^{\qQ}$
 implies the functor $T_{\eE}\colon D(X)\to D(X)$
 with kernel
$\rR$ fits into the triangle
  $$\dR \Hom _X(\pi _{\ast}\eE, F) \dotimes _R \pi _{\ast}\eE \lr 
F \lr T_{\eE}(F), $$ 
for $F\in D(X)$. 
We check $T_{\eE}$ gives an equivalence. 
We follow the arguments of~\cite{Pl},~\cite{HT}. 
Define $E^{\perp}$ to be the subcategory
 $\{ F\in D(X) \mid \dR \Hom (E, F)=0 \}$. 
Then $\Omega \cneq 
E\cup E^{\perp}$ is a spanning class. 
Let $<E>$ be the minimum extension closed subcategory of $D(X)$ which 
contains $E$. Then since $R$ is finite dimensional,
 we have $\pi _{\ast}\eE 
\in <E>$. 
Therefore if $F\in E^{\perp}$, then 
$\dR \Hom _X(\pi _{\ast}\eE, F)=0$. 
Hence $T_{\eE}(F)\cong F$ for $F\in E^{\perp}$. Next since 
$\eE$ is $R$-spherical, we have the distinguished triangle:
$$E \lr \dR \Hom _X (\pi _{\ast}\eE, E)\dotimes _R \pi _{\ast}\eE 
\lr E[-\dim X]. $$
Then the following diagram
$$\xymatrix{
E \ar[d]\ar[dr]^{\id} & & \\
\dR \Hom _X (\pi _{\ast}\eE, E)\dotimes _R \pi _{\ast}
\eE \ar[d] \ar[r] & E \ar[d]
\ar[r] & 
T_{\eE}(E) \\
E[-\dim X] \ar[r] & 0 \ar[ur] & \\
}$$
shows $T_{\eE}(E)\cong E[1-\dim X]$. 
Therefore $T_{\eE}$ is fully faithful on $\Omega$, 
hence fully faithful on $D(X)$. (\textit{c.f.}~\cite[Theorem 2.3]{Br2}).
Finally the assumption $E\otimes \omega _X \cong E$ implies 
$T_{\eE}(F)\otimes \omega _X \cong T_{\eE}(F)$ for $F\in \Omega$. 
Therefore $T_{\eE}$ gives an equivalence by the 
argument of~\cite[Theorem 5.4]{Br2}. 
$\quad \square$

\section{Flops at $(0, -2)$-curves}
We give some examples of autoequivalences associated to $R$-spherical 
objects. 
Let $f\colon X\to Y$ be a three dimensional flopping 
contraction which contracts a rational curve $C\subset X$.
Let $f^{\dag}\colon X^{\dag}\to Y$ be its flop, 
and $C^{\dag}\subset X^{\dag}$ be the flopped curve. 
Then in~\cite{B-O2},~\cite{Br1},~\cite{Ch}, the functor
 $\Phi _1 \colon D(X^{\dag})\to D(X)$
with kernel $\oO _{X\times _Y X^{\dag}}$ gives an equivalence.  
 $\Phi _1$ satisfies the 
following (\textit{c.f.}\cite[Lemma 5.1]{Tst}):
\begin{itemize}
\item $\Phi _1$ takes $\oO _{C^{\dag}}(-1)[1]$ to $\oO _C (-1)$. 
\item $\Phi _1$ commutes with derived push-forwards. i.e. 
$\dR f_{\ast}\circ \Phi _1 \cong \dR f^{\dag}_{\ast}$. 
\end{itemize}
Similarly we can construct the equivalence 
$\Phi _2 \colon D(X)\to D(X^{\dag})$ with kernel 
$\oO _{X\times _Y X^{\dag}}$. 
Composing these, we
obtain the autoequivalence 
$$\Phi \cneq \Phi _1 \circ \Phi _2 
\colon D(X)\lr D(X^{\dag}) \lr  D(X).$$
 Note that $\Phi (\oO _C (-1))=\oO _C (-1)[-2]$ and 
$\Phi$ commutes with $\dR f_{\ast}$.  If $C\subset X$ is a $(-1, -1)$-curve, 
then $\oO _C (-1)$ is a spherical object and $\Phi$ coincides with 
the associated twist $T_{\oO _C (-1)}$. But if $C$ is not a $(-1, -1)$-curve, 
then $\oO _C (-1)$ is no longer spherical, so we have to find some 
new descriptions of $\Phi$. The idea is to consider the moduli 
problem of $\oO _C (-1)$ and using the universal family. 

Here we 
assume $C\subset X$ is a $(0, -2)$-curve, i.e. 
normal bundle is $\oO _{C}\oplus \oO _C (-2)$,  and 
give the description of $\Phi$. 
Let $\mM$ be the connected component of the moduli space of simple 
sheaves on $X$, which contains $\oO _C (-1)$. 
We define $R_m$, $S_m$ to be 
$$R_m \cneq \mathbb{C}[t]/(t^{m+1}), \qquad S_m \cneq \Spec R_m .$$
Since $\Ext _X ^1 (\oO _C (-1), \oO _C (-1))=\mathbb{C}$ and 
$C\subset X$ is rigid, we can write $\mM$ as 
$\mM = S_m$ for some $m\in \mathbb{N}$. 
Let $\eE \in \Coh (X\times S_m)$ be the universal family. 

\begin{thm}\label{flop}
$\eE$ is a $R_m$-spherical object and the associated functor $T_{\eE}\colon 
D(X)\to D(X)$ coincides with $\Phi$. 
\end{thm}
\textit{Proof}. 
For $n\le m$, define $\eE _n$ to be 
$$\eE _n \cneq 
\pi _{n\ast}(\eE |_{X\times S_n}) \in \Coh (X),$$
where $\pi _{n}\colon X\times S_n \to X$ is a projection. 
Since we have the exact sequences of $R_m$-modules:
\begin{align*}
& 0 \lr R_{n-1} \lr R_n \lr \mathbb{C} \lr 0, \\
& 0 \lr \mathbb{C}\lr R_n \lr R_{n-1}\lr 0, 
\end{align*}
we have the exact sequences in $\Coh (X)$:
\begin{align}
& 0 \lr \eE_{n-1} \lr \eE_n \lr E \lr 0, \\
& 0 \lr E \lr \eE_{n} \lr \eE_{n-1}\lr 0.
\end{align}
Applying $\Hom (\ast, \eE)$ to the sequence $(1)$, we obtain the 
long exact sequence:
\begin{align*}
\Hom (\eE _n , E) \lr \Hom (\eE _{n-1}, E)  \stackrel{\xi _n}{\lr}
\Ext ^1(E, E)=\mathbb{C} \\
\lr \Ext ^1 (\eE _n, E) \lr \Ext ^1 (\eE _{n-1}, E) \stackrel{\eta _n}{\lr}
\Ext ^2 (E, E)=\mathbb{C}.
\end{align*}
On the other hand, the
 sequence $(2)$ determines the non-zero element:
$$e_n \in \Ext ^1 (\eE _{n-1}, E), $$
and $\eta _n (e_n)\in \Ext ^2 (E, E)$ gives the obstruction 
to deforming $\eE |_{X\times S_n}$ to a coherent sheaf on 
$X\times S_{n+1}$ flat over $S_{n+1}$. 
(\textit{c.f.}~\cite[Proposition 3.13]{Thom}).
Therefore $\eta _n (e_n)=0$ for $n<m$ and $\eta _m (e_m)\neq 0$. 
On the other hand, we have the following morphism of exact sequences:
$$\xymatrix{
0 \ar[r] & \eE _{n-1} \ar[r] \ar[d]_{s_n}
 & \eE _n \ar[r]\ar[d] & E \ar[r]\ar[d] & 
0 \\
0 \ar[r] & E \ar[r] & \eE _1 \ar[r] & E \ar[r] & 0,
}$$
where $s_n$ is a natural surjection. 
Hence $\xi _n (s _n) \in \Ext ^1 (E, E)$ corresponds to the 
extension $\eE _1$, which is a non-trivial first order deformation of 
$E$. Therefore $\xi _n (s _n)\neq 0$ and $\xi _n$ is surjective. 
Combining these, we have 
\begin{align*}
\Ext ^1 (\eE _{n-1}, E) &\cong \Ext ^1 (\eE _{n}, E) \\
& \cong \mathbb{C} \quad (\mbox{ for } n<m), \\
\Ext ^1 (\eE _m, E) &\cong 0, \\
\Hom (\eE _{n-1}, E) & \cong \Hom (\eE _n, E) \cong \mathbb{C}. 
\end{align*}
Similarly applying $\Hom (E, \ast)$ to the sequence $(2)$,
we obtain $\Ext ^1 (E, \eE _m)=0$ and $\Hom (E, \eE _m)=\mathbb{C}$. 
By Serre duality, we can conclude $\eE$ is $R_m$-spherical.  

Next let us consider the equivalence:
$$\widetilde{\Phi}  \cneq T_{\eE}\circ \Phi ^{-1}\colon D(X) \lr D(X).$$
Then $\widetilde{\Phi}$ takes $\oO _C (-1)$ to $\oO _C (-1)$, 
and commutes with $\dR f_{\ast}$. Therefore $\widetilde{\Phi}$ 
preserves perverse t-structure $\oPPer (X/Y)$ in the 
sense of~\cite{Br1}. 
Then the argument of~\cite[Theorem 6.1]{Tst}
shows $\widetilde{\Phi}$ is isomorphic to the identity functor.
 $\quad \square$
 
 \section{Deformations of $\mathbb{P}$-twists}
 \subsection*{Review of $\mathbb{P}$-objects and associated twists}
 $R$-spherical twists can also be used to 
 construct deformations of $\mathbb{P}$-twists.
Let us recall the definition of $\mathbb{P}$-objects and 
the associated autoequivalences introduced in~\cite{HT}. 
Again we assume $X$ is a smooth projective variety over $\mathbb{C}$. 
\begin{defi}\emph{\bf{\cite{HT}}}
An object $E\in D(X)$ is called $\mathbb{P}^n$-object 
if it satisfies the following:
\begin{itemize}
\item $\Ext ^{\ast}_X (E, E)$ is isomorphic to 
$H^{\ast}(\mathbb{P}^n, \mathbb{C})$ as a graded ring.   
\item $E\otimes \omega _X \cong E$. 
\end{itemize}
\end{defi}
Note that if $\mathbb{P}^n$-object exists, then $\dim X =2n$ 
by Serre duality. D.Huybrecht and R.Thomas~\cite{HT}
 constructed an 
equivalence $P_E \colon D(X)\to D(X)$ associated to $E$, which 
is described as follows. 
Let $h\in \Ext _X ^2 (E, E)$ be the degree two generator. 
First consider the 
morphism in $D(X\times X)$:
$$H\cneq \check{h}\boxtimes \id - \id \boxtimes h \colon 
\check{E}\boxtimes E[-2] {\lr}
\check{E}\boxtimes E. $$
Let us take its cone $\hH \in D(X\times X)$. We can see
 the composition $H$ with 
the trace map $\tr \colon \check{E}\boxtimes E\to \oO _{\Delta}$ 
becomes zero. Therefore there exists a (in fact unique) morphism
$t \colon \hH \to \oO _{\Delta}$ such that the following diagram commutes~\cite[Lemma 2.1]{HT}:
$$\xymatrix{
\check{E}\boxtimes E [-2] \ar[r]^{H} & \check{E}\boxtimes E \ar[r]\ar[d]_{\tr}
 & \hH \ar[ld]^{t} \\
 & \oO _{\Delta}. & 
 }$$
 Then define $\qQ _{\eE}$ to be the cone
 $$\qQ _{\eE}\cneq \Cone (t\colon \hH \lr \oO _{\Delta}) \in D(X\times X).$$
 Then in~\cite{HT}, it is shown that 
  the functor $P_{\eE}\colon D(X)\to D(X)$ with kernel 
 $\qQ _{\eE}$ gives the equivalence. 
 
 Next let us consider a one parameter deformation of $X$. 
 Let $f\colon \mathrsfs{X}\to C$ be a smooth family over a 
 smooth curve $C$ with a distinguished fibre
  $j\colon X= f^{-1}(0)\hookrightarrow \mathrsfs{X}$,
  $0\in C$. 
Suppose $E\in D(X)$ is a $\mathbb{P}^n$-object and let
$A(E)\in \Ext _X ^1(E, E\otimes \Omega _X)$ be its Atiyah-class. 
Then the obstruction to deforming $E$ sideways to 
first order is given by the product
$$A(E)\cdot \kappa (X)\in \Ext _X ^2 (E, E), $$
where $\kappa (X)\in H^1 (X, T_X)$ is the Kodaira-Spencer class of 
the family $f\colon \mathrsfs{X} \to C$. In~\cite{HT},
the case of $A(E)\cdot \kappa (X)\neq 0$ is studied. 
In that case, $j_{\ast}E$ is a spherical object and the 
associated equivalence $T_{j_{\ast}E}\colon D(\mathrsfs{X})\to D(\mathrsfs{X})$ fits into the 
commutative diagram~\cite[Proposition 2.7]{HT}
$$\xymatrix{
D(X) \ar[r]^{j_{\ast}}\ar[d]_{P_E} & D(\mathrsfs{X}) \ar[d]^{T_{j_{\ast}E}} \\
D(X) \ar[r]_{j_{\ast}} & D(\mathrsfs{X}). 
}$$
Our purpose is to treat the case of $A(E)\cdot \kappa (X)= 0$.

\subsection*{$R$-spherical objects via deformations of $\mathbb{P}$-objects} 
Let $f\colon \mathrsfs{X} \to C$ and $E\in D(X)$ be as before, and 
assume $A(E)\cdot \kappa (X)=0$. 
Note that $j_{\ast}E$ is not spherical. 
In fact we have the 
distinguished triangle
$$E[1] \lr \dL j^{\ast}j_{\ast}E \lr E \stackrel{A(E)\cdot \kappa (X)}{\lr}
E[2],$$
by~\cite[Proposition 3.1]{HT}. Hence we have the decomposition
$\dL j^{\ast}j_{\ast}E \cong E\oplus E[1]$, and we calculate 
\begin{align*}
\Ext _{\mathrsfs{X}}^k (j_{\ast}E, j_{\ast}E) & \cong \Ext _X ^k (\dL j^{\ast}j_{\ast}E, E)\\
& \cong \Ext _X ^k (E, E)\oplus \Ext _X ^{k-1}(E, E) \\
& \cong \mathbb{C}
\end{align*}
for $0\le k\le 2n+1$. 
As in the previous section, we are going to consider deformations of 
$j_{\ast}E$ in $\mathrsfs{X}$. 
The moduli theories of complexes were carried out 
by~\cite{Inaba},~\cite{LIE}. 
Following the notation used in~\cite{Inaba}, we consider the functor
$\Splcpx _{\mathrsfs{X} /C}$ from the category of locally 
noetherian schemes over 
$C$ to the category of sets, 
$$\Splcpx _{\mathrsfs{X} /C}(T)\cneq \left\{ 
\fF ^{\bullet} \ \begin{array}{|l} \fF ^{\bullet}
\mbox{ is a bounded complex of coherent sheaves on }\mathrsfs{X} _T
 \\ \mbox{ such that each }
\fF ^i \mbox{ is flat over }T \mbox{ and for any }
t\in T, \\ \Ext _{X_t}^0 (\fF ^{\bullet}(t), \fF ^{\bullet}(t))\cong k(t), 
\Ext _{X_t}^{-1} (\fF ^{\bullet}(t), \fF ^{\bullet}(t))=0 
\end{array} \right\}/\sim .$$
Here $\mathrsfs{X} _T \cneq \mathrsfs{X} \times _{C}T$, $\fF ^{\bullet}(t)\cneq \fF ^{\bullet}
\otimes _T k(t)$, and $\fF ^{\bullet}\sim \fF ^{'\bullet}$ if and only if
there exist $\lL \in \Pic (T)$, a bounded complex of quasi-coherent 
sheaves $\gG ^{\bullet}$ and quasi-isomorphisms 
$\gG ^{\bullet}\to \fF ^{\bullet}$, $\gG ^{\bullet}\to \fF ^{'\bullet}\otimes 
\lL$. Let $\Splcpx _{\mathrsfs{X} /C}^{\rm{et}}$ be the associated sheaf
of $\Splcpx _{\mathrsfs{X} /C}$ in the \'{e}tale topology. M.Inaba~\cite{Inaba}
showed the following:
\begin{thm}\emph{\bf{\cite{Inaba}}}
The functor
$\Splcpx _{\mathrsfs{X} /C}^{\emph{et}}$ is represented by a locally separated 
algebraic space $\mM$ over $C$. 
\end{thm}
Let $\gamma \colon S_1 \hookrightarrow C$ 
be an extension of $0\hookrightarrow C$. Let $r$ be the 
restriction, 
$$r\colon \Splcpx _{\mathrsfs{X} /C}^{\emph{et}}(\gamma) \lr \Splcpx _{\mathrsfs{X} /C}^{\emph{et}}(0).$$
By the assumption $A(E)\cdot \kappa (X)=0$, we have 
$r^{-1}(E)\neq \emptyset$. Moreover by~\cite[Proposition 2.3]{Inaba}, 
there is a bijection between $r^{-1}(E)$ and $\Ext _X ^1 (E, E)$, which 
is zero. 
Therefore the map $T_{\mM, E} \to T_{C, 0}$ is an isomorphism, hence
$\dim \mM \le 1$ at $[E]\in \mM$. 
Note that by taking push-forward along the inclusion
$\mathrsfs{X}\times _C T \to \mathrsfs{X}
\times T$, 
we get the morphism of functors:
$$\delta \colon 
\Splcpx _{\mathrsfs{X} /C}^{\emph{et}} \lr \Splcpx _{\mathrsfs{X} /S_0}^{\emph{et}}.$$
We put the following technical assumption $(\star)$. 
\begin{itemize}
\item The morphism $\delta$ gives an isomorphism between connected components 
of both sides, 
which contain $E$ and $j_{\ast}E$ respectively.
Let $[E]\in \mM' \subset \mM$ be the connected component. We assume 
$\mM '$ is a zero-dimensional scheme.  $\cdots (\star)$. 
\end{itemize}

Note that we can write $\mM '=S_m$ for 
some $m$.  
Let $\mathrsfs{X} _m 
\cneq \mathrsfs{X} \times _C \mM ' =\mathrsfs{X}\times _C S_m$ and 
$\eE \in D(\mathrsfs{X} _m)$ be the universal family. 
We use the following notations for morphism:
$$
\xymatrix{
\mathrsfs{X} _m \ar[d]_{f'} \ar[r]^{k} & \mathrsfs{X} \ar[d]^{f} \\
S_m \ar[r]_{k'} & C ,
} \quad
\xymatrix{
\mathrsfs{X} _m \ar[r]^{l} & \mathrsfs{X} \times S_m \ar[rd]_{\pi} 
\ar[r]^{\pi '} & S_m  \\
X \ar[u]_{i} \ar[r]_{j} & \mathrsfs{X} \ar[u]_{i'} & \mathrsfs{X}.
}
$$
If there is no confusion, we will use the same notations for 
$n\le m$. 
We show the following proposition:

\begin{prop}\label{spherical}
The object
$l_{\ast}\eE \in D(\mathrsfs{X} \times S_m)$ is $R_m$-spherical. 
\end{prop}
\textit{Proof}. 
Since $\pi _{\ast}l_{\ast}\eE \cong k_{\ast}\eE$
and $\dL i^{'\ast}l_{\ast}\eE \cong j_{\ast}E$, 
we have to calculate $\Ext _{\mathrsfs{X}}^i (k_{\ast}\eE, j_{\ast}E)$. 
By the assumption $(\star)$, 
we cannot deform $l_{\ast}\eE$ to $(m+1)$-th order. For $n\le m$, let
$\eE _n \cneq \eE |_{\xX _n} \in D(\xX _n)$ and
$\widetilde{\eE}_n \cneq k_{\ast}\eE _n \in D(\mathrsfs{X})$. We 
consider distinguished triangles:
\begin{align}
& \widetilde{\eE}_{n-1} \lr \widetilde{\eE}_n \lr j_{\ast}E \stackrel{e_n '}
\lr \widetilde{\eE}_{n-1}[1], \\
& j_{\ast}E
 \lr \widetilde{\eE}_{n} \lr \widetilde{\eE}_{n-1}\stackrel{e_n}{\lr}
j_{\ast}E[1].
\end{align}
Then by the argument of~\cite[Proposition 3.3]{Thom}, 
we can see that
the 
composition
$$e_n \circ e_n ' \colon j_{\ast}E \lr \widetilde{\eE}_{n-1}[1] \lr 
j_{\ast}E[2]$$
gives the obstruction to deforming $l_{\ast}\eE _n$ to $(n+1)$-th 
order. If $E$ is a sheaf, this is just~\cite[Proposition 3.3]{Thom}
and we can generalize this by replacing
 the exact sequences in~\cite[Proposition 3.3]{Thom}
by the exact sequences of representing complexes. We leave the detail 
to the reader. 
Hence $e_m \circ e_m ' \neq 0$ and $e_n \circ e_n ' =0$ for 
$n<m$. Applying $\Hom (\ast, j_{\ast}E)$ to the triangle $(3)$, 
we obtain the long exact sequence, 
$$\Ext _{\mathrsfs{X}}^1 (j_{\ast}E, j_{\ast}E) \lr \Ext _{\mathrsfs{X}}^1 (\widetilde{\eE}_n, 
j_{\ast}E) \lr \Ext _{\mathrsfs{X}}^1 (\widetilde{\eE}_{n-1}, j_{\ast}E) 
\lr \Ext _{\mathrsfs{X}}^2 (j_{\ast}E, j_{\ast}E)=\mathbb{C}. $$
Then using the above sequence and the same argument as in Theorem~\ref{flop}, 
we can conclude $\Ext _{\mathrsfs{X}}^1 (k_{\ast}\eE, j_{\ast}E)=0$. 

Next we use the existence of the distinguished triangle 
below~\cite[Lemma 3.3]{B-O2}:
$$\eE [1] \lr \dL k^{\ast}k_{\ast}\eE \lr \eE \lr \eE [2].$$
Pulling back to $X$, we have the triangle:
\begin{align}
E[1] \lr \dL j^{\ast}k_{\ast}\eE \lr E \stackrel{\theta}{\lr} E[2]. \end{align}
Since $\Ext _X ^2 (E,E)$ is one dimensional, $\theta$ is zero or non-zero
multiple of $h$. Assume $\theta =0$. Then $\dL j^{\ast}k_{\ast}\eE
\cong E \oplus E[1]$, and 
\begin{align*}
\Ext _{\mathrsfs{X}}^1 (k_{\ast}\eE, j_{\ast}E) & \cong \Ext _{\mathrsfs{X}}^1 (\dL j^{\ast}k_{\ast}\eE, E) \\
& \cong \Ext _X ^1 (E, E) \oplus \Hom (E, E) \\
& \cong \mathbb{C}, 
\end{align*}
which is a contradiction. Hence we may assume $\theta =h$. 
Applying $\Hom (\ast, E)$ to the triangle $(5)$, 
we obtain the long exact sequence:
$$\lr \Ext _X ^i (E, E) \lr \Ext _{\mathrsfs{X}}^i (\dL j^{\ast}k_{\ast}\eE, E) \lr
\Ext _X ^{i-1}(E, E) \stackrel{h}{\lr} \Ext _X ^{i+1}(E, E)\lr .$$
By the definition of $\mathbb{P}^n$-object, we obtain 
\begin{align*}
\Ext _{\mathrsfs{X}}^i (k_{\ast}\eE, j_{\ast}E) & \cong \Ext _X ^i (\dL j^{\ast}k_{\ast}
\eE, E) \\
& =\left\{ \begin{array}{cl} \mathbb{C} &(i=0, 2n+1), \\
0 & (i\neq 0, 2n+1). \quad \square \end{array} \right. 
\end{align*}

\vspace{3mm}

\begin{rmk}
The assumption $(\star)$ is satisfied if $E$ is a sheaf
and $\dim \mM ' =0$. In fact suppose $l_{\ast}\eE$ 
extends to a $S_{m+1}$-valued point of 
$\Splcpx ^{\emph{et}}_{\mathrsfs{X}/S_0}$.  
Then as in\emph{~\cite[Proposition 3.13]{Thom}},
there exists $\widetilde{\eE}_{m+1}\in \Coh (\mathrsfs{X})$
such that there exists a morphism of exact sequences of
 $\oO _{\mathrsfs{X}}$-modules:
$$\xymatrix{
0 \ar[r] & \widetilde{\eE}_{m-1} \ar[r] & \widetilde{\eE}_{m} \ar[r] & 
j_{\ast}E \ar[r] & 0 \\
0 \ar[r]  & \widetilde{\eE}_{m} \ar[r]^{\nu}\ar[u]
 & \widetilde{\eE}_{m+1} \ar[r]\ar[u] & 
j_{\ast}E \ar[r]\ar[u]_{\emph{id}} & 0.
}
$$
An easy diagram chasing 
shows $\widetilde{\eE}_{m+1}$ is a $\oO _{\mathrsfs{X}}/(t^{m+2})$-module
for the uniformizing parameter $t\in \oO _{C,0}$. 
Moreover we have $t\cdot \widetilde{\eE}_{m+1}=\Imm \nu$. Therefore 
the map
$$\widetilde{\eE}_{m+1}\otimes _{\oO _C/(t^{m+2})}(t) \lr
 \widetilde{\eE}_{m+1}$$
 is a morphism from $\widetilde{\eE}_{m}$ 
 onto $\Imm \nu \cong \widetilde{\eE}_m$, hence injective. 
 Then 
\emph{~\cite[Lemma 3.7]{Thom}} shows
$\widetilde{\eE}_{m+1}$ 
is flat over $\oO _{C,0}/(t^{m+2})$ and gives a $S_{m+1}$-valued
point of $\Splcpx ^{\emph{et}}_{\mathrsfs{X}/C}$.
\end{rmk}

\subsection*{$\mathbb{P}$-twists and $R$-spherical twists}
By Proposition~\ref{spherical}, we have the associated functor 
$T_{l_{\ast}\eE}\colon D(\mathrsfs{X})\to D(\mathrsfs{X})$ 
under the assumption $(\star)$. The next 
purpose is 
to show the existence of the diagram as in~\cite[Proposition 2.7]{HT}.  
We use the following notations for morphisms:
$$
\xymatrix{
X\times X \ar[dr]_{\widetilde{i}} \ar[rr]^{\widetilde{j}} &  &  
\mathrsfs{X} \times _C \mathrsfs{X}  \\
& \mathrsfs{X} _m \times _{S_m}\mathrsfs{X} _m , \ar[ur]_{\widetilde{k}} &
}
\xymatrix{
\mathrsfs{X} \times \mathrsfs{X} _m \ar[r]^{\id \times l}&  \mathrsfs{X} \times S_m \times \mathrsfs{X} \\
\mathrsfs{X} _m \times _{S_m}\mathrsfs{X} _m \ar[u]^{l'} \ar[r]_{l''} & \mathrsfs{X} _m \times \mathrsfs{X} ,
\ar[u]_{l\times \id}
}
$$

$$
\xymatrix{
X\times X \ar[r]^{\widetilde{i}} & \mathrsfs{X} _m \times _{S_m} \mathrsfs{X} _m \ar[r]^{\widetilde{k}} & \mathrsfs{X} \times _C \mathrsfs{X} \ar[r]^{\iota} & \mathrsfs{X} \times \mathrsfs{X} \\
X \ar[u]_{\Delta _0} \ar[r]_{i} & \mathrsfs{X} _m  \ar[u]_{\Delta _m} \ar[r]_{k} & 
\mathrsfs{X} , \ar[u]_{\Delta '} \ar[ur]_{\Delta } 
}$$
$$
\xymatrix{
& \mathrsfs{X} \times S_m \times 
\mathrsfs{X} \ar[dl]_{p_{12}} \ar[d]_{p_{13}} \ar[dr]^{p_{23}} & \\
\mathrsfs{X} \times S_m & \mathrsfs{X} \times \mathrsfs{X} & \mathrsfs{X} \times S_m , 
} \quad 
\xymatrix{
& \mathrsfs{X} _m \times _{S_m} \mathrsfs{X} _m \ar[dl]_{q_1} \ar[dr]^{q_2} & \\
\mathrsfs{X} _m & & \mathrsfs{X} _m .
}$$

\begin{thm}\label{com}
The functor
$T_{l_{\ast}\eE}$ fits into the following commutative diagram:
$$\xymatrix{
D(X) \ar[r]^{j_{\ast}} \ar[d]_{P_{E}} & D(\mathrsfs{X})\ar[d]^{T_{l_{\ast}\eE}} \\
D(X) \ar[r]_{j_{\ast}} & D(\mathrsfs{X}). }$$
\end{thm}
\textit{Proof}. 
We try to 
imitate the argument of~\cite[Proposition 2.7]{HT}. 
First we construct the morphism 
$$\alpha \colon \widetilde{k}_{\ast}(q_1 ^{\ast}\check{\eE} \dotimes 
q_2 ^{\ast}\eE )[-1]\lr \Delta ' _{\ast} \oO _{\mathrsfs{X}}$$
in $D(\mathrsfs{X} \times _C \mathrsfs{X})$. This is constructed 
 by the composition of 
$\widetilde{k}_{\ast}\tr$, 
$$\widetilde{k}_{\ast}\tr \colon \widetilde{k}_{\ast}
(q_1 ^{\ast}\check{\eE}\dotimes q_2 ^{\ast}\eE)[-1] \lr \widetilde{k}_{\ast}
\Delta _{m\ast}\oO _{\mathrsfs{X} _m}[-1]=\Delta ' _{\ast}k_{\ast}\oO _{\mathrsfs{X} _m}[-1], $$
with the morphism $\Delta ' _{\ast}k_{\ast}\oO _{\mathrsfs{X} _m}[-1] \to 
\Delta ' _{\ast}
\oO _{\mathrsfs{X}}$ obtained 
 by applying $\Delta ' _{\ast}$ to the exact sequence, 
$$0 \lr \oO _{\mathrsfs{X}} \lr \oO _{\mathrsfs{X}}(\mathrsfs{X} _m) \lr k_{\ast}\oO _{\mathrsfs{X} _m}\lr 0.$$
Let $\lL \cneq \Cone (\alpha) \in D(\mathrsfs{X} \times _C \mathrsfs{X})$. 
Applying Chen's lemma~\cite{Ch}, it suffices to show 
\begin{align*}
\iota _{\ast}\lL &\cong
 \Cone \left( \dR p_{13\ast}(p_{12}^{\ast}
(\check{l_{\ast}\eE}
\dotimes \pi ^{!}\oO _{\mathrsfs{X}}) \dotimes p_{23}^{\ast}
l_{\ast}\eE) \stackrel{\mu}{\lr} \Delta _{\ast}\oO _{\mathrsfs{X}} \right), \\
\dL \widetilde{j}^{\ast}\lL &\cong \hH .
\end{align*}
Here $\mu$ is the morphism constructed in the proof of Theorem~\ref{equ}
and $\hH$ is the kernel of $P_{E}$. First we check $\iota _{\ast}\lL
\cong \Cone (\mu)$. Note that $\pi ^{!}\oO _{\mathrsfs{X}}=\oO _{\mathrsfs{X} \times S_m}$
and $\check{l_{\ast}\eE} \cong l_{\ast}\check{\eE}[-1]$ by the 
duality isomorphism. Hence we have 
\begin{align*}
\dR p_{13\ast}(p_{12}^{\ast}
(\check{l_{\ast}\eE}\dotimes \pi ^{!}\oO _{\mathrsfs{X}}) \dotimes p_{23}^{\ast}
l_{\ast}\eE) 
& \cong \dR p_{13\ast}(p_{12}^{\ast}l_{\ast}\check{\eE} \dotimes 
p_{23}^{\ast}l_{\ast}\eE)[-1] \\
& \cong \dR p_{13\ast}\left\{ (l\times \id)_{\ast}r_1 ^{\ast}\check{\eE}
\dotimes (\id \times l)_{\ast}r_2 ^{\ast}\eE \right\}[-1] \\
& \cong \dR p_{13\ast}(\id \times l)_{\ast}\left\{
\dL (\id \times l)^{\ast}(l\times \id)_{\ast}r_1 ^{\ast}\check{\eE}\dotimes 
r_2 ^{\ast}\eE \right\}[-1] \\
& \cong \dR p_{13\ast}(\id \times l)_{\ast}(l' _{\ast}\dL l^{''\ast}r_1 ^{\ast}
\check{\eE}\dotimes r_2 ^{\ast}\eE )[-1] \\
& \cong \dR p_{13\ast}(\id \times l)_{\ast}l' _{\ast}(\dL l^{''\ast} 
r_1 ^{\ast}\check{\eE}\dotimes \dL l^{'\ast}r_2 ^{\ast}\eE )[-1] \\
& \cong \iota _{\ast}\widetilde{k}_{\ast}(q_1 ^{\ast}\check{\eE}
\dotimes q_2 ^{\ast}\eE)[-1]. 
\end{align*}
Here $r_1$, $r_2$ are defined by the fiber squares:
$$\xymatrix{
\mathrsfs{X} _m \times \mathrsfs{X} \ar[r]^{l\times \id} \ar[d]_{r_1} & \mathrsfs{X} \times S_m \times 
\mathrsfs{X} \ar[d]_{p_{12}} \\
\mathrsfs{X} _m \ar[r] _{l} & \mathrsfs{X} \times S_m ,
} \quad 
\xymatrix{
\mathrsfs{X} \times \mathrsfs{X} _m \ar[r]^{\id \times l} \ar[d]_{r_2} & \mathrsfs{X} \times S_m \times 
\mathrsfs{X} \ar[d]_{p_{23}} \\
\mathrsfs{X} _m \ar[r] _{l} & \mathrsfs{X} \times S_m .
}$$

Under the above isomorphism, we can check 
$\iota _{\ast}\alpha =\mu$. 
Hence we have 
$\widetilde{l}_{\ast}\lL \cong \Cone (\mu)$. 

Next we check $\dL \widetilde{j}^{\ast}\lL \cong \hH$. 
Note that we have 
$$\dL \widetilde{j}^{\ast}\lL 
= \Cone \left( \dL \widetilde{j}^{\ast}
\widetilde{k}_{\ast}(q_1 ^{\ast}\check{\eE}\dotimes q_2 ^{\ast}\eE )[-1]
\stackrel{\dL \widetilde{j}^{\ast}\alpha}{\lr} \dL \widetilde{j}^{\ast}
\Delta _{\ast}' \oO _{\mathrsfs{X}}=\Delta _{0\ast}\oO _{X} \right), $$
and there exists the distinguished triangle
$$q_1 ^{\ast}\check{\eE}\dotimes q_2 ^{\ast}\eE [-2] \lr 
q_1 ^{\ast}\check{\eE}\dotimes q_2 ^{\ast}\eE 
\lr \dL \widetilde{k}^{\ast}\widetilde{k}_{\ast}(q_1 ^{\ast}\check{\eE}
\dotimes q_2 ^{\ast}\eE )[-1],$$
as in~\cite[Lemma 3.3]{B-O2}. Then applying $\dL \widetilde{i}^{\ast}$, 
we have the triangle:
$$\check{E}\boxtimes E [-2] \stackrel{u}{\lr} 
\check{E}\boxtimes E \lr 
\dL \widetilde{j}^{\ast}\widetilde{k}_{\ast}(q_1 ^{\ast}\check{\eE}
\dotimes q_2 ^{\ast}\eE)[-1]. $$
We can easily check the following:
\begin{align*}
& \Ext _{X\times X}^2 (\check{E}\boxtimes E , \check{E}\boxtimes E) \\
& \cong \left( \Ext _X ^2 (E, E)\otimes \Ext _X ^0 (E, E) \right)\oplus 
\left(
\Ext _X ^0 (E, E)\otimes \Ext _X ^2 (E, E) \right) . 
\end{align*}
Hence we can write $u=a(\check{h}\boxtimes \id)+ b(\id \boxtimes h)$
for some $a, b\in \mathbb{C}$. 
On the other hand, we can check that the diagram
$$\xymatrix{
\check{E}\boxtimes E [-2] \ar[r]^{u} &
\check{E}\boxtimes E \ar[r]\ar[dr]_{\tr}
 & \dL \widetilde{j}^{\ast}\widetilde{k}_{\ast}
(q_1 ^{\ast}\check{\eE}\dotimes q_2 ^{\ast}\eE) 
\ar[d] ^{\dL \widetilde{j}^{\ast}\alpha} \\
& & \Delta _{0\ast}\oO _X, 
 }$$
 commutes. 
 This is easily checked using the
  same argument of~\cite[Proposition 2.7]{HT}, and leave the detail
  to the reader. 
 Therefore $\tr \circ u=0$, which implies $b=-a$. 
 Hence if we show $u\neq 0$, then we can conclude $\dL \widetilde{j}^{\ast}
 \lL \cong \hH$. 
 Assume $u=0$. Then we have the decomposition
 \begin{align}
 \dL \widetilde{j}^{\ast}\widetilde{k}_{\ast}(q_1 ^{\ast}\check{\eE}
 \dotimes q_2 ^{\ast}\eE)[-1] \cong (\check{E}\boxtimes E )\oplus 
( \check{E}\boxtimes E)[-1]. \end{align}
 Since 
 $\Hom _{X\times X}(\check{E}\boxtimes E [-1], \Delta _{\ast}\oO _X)=0,$
 the morphism 
 $$\dL \widetilde{j}^{\ast}\alpha \colon 
 \dL \widetilde{j}^{\ast}\widetilde{k}_{\ast}(q_1 ^{\ast}\check{\eE}\dotimes 
 q_2 ^{\ast}\eE)[-1] \lr \Delta _{0\ast}\oO _X$$
 is a non-zero multiple of $(\tr, 0)$ under the decomposition $(6)$.  
 Let $\mathrsfs{S} \in D(X\times X)$ be the cone of the trace map:
 $$\check{E}\boxtimes E \stackrel{\tr}{\lr} \Delta _{0\ast}\oO _X 
 \lr \mathrsfs{S}.$$
 Then we have the decomposition
 $\dL \widetilde{j}^{\ast}\lL 
 \cong \mathrsfs{S} \oplus (\check{E}\boxtimes E)$,
 and the following diagram commutes:
 $$\xymatrix{
    D(X) \ar[r]^{\Phi _{X\to X}^{\dL \widetilde{j}^{\ast}\lL}} 
 \ar[d]_{j\ast} & 
 D(X) \ar[d]^{j\ast} \\
  D(\mathrsfs{X}) \ar[r]_{T_{l_{\ast}\eE}} & D(\mathrsfs{X}), }$$
  by Chen's lemma~\cite{Ch}.
In particular we have 
\begin{align*}
j_{\ast}\Phi _{X\to X}^{\dL \widetilde{j}^{\ast}\lL}(E) & \cong T _{l_{\ast}\eE}
(j_{\ast}E) \\
& \cong j_{\ast}E [1-\dim \mathrsfs{X}], \end{align*}
which is indecomposable. It follows that 
$$\Phi _{X\to X}^{\mathrsfs{S}}(E)\cong 
0 \quad \mbox{or} \quad \Phi _{X\to X}^{\check{E}
\boxtimes E}(E) \cong 0.$$
Since $\Phi _{X\to X}^{\check{E}
\boxtimes E}(E) \cong \dR \Hom (E, E)\otimes _{\mathbb{C}}E$, the latter 
is impossible by the definition of $\mathbb{P}^n$-object.
 Hence $\Phi _{X\to X}^{\mathrsfs{S}}(E)$ must be zero. Since we have the distinguished 
 triangle:
 $$\dR \Hom (E, E)\otimes _{\mathbb{C}}E \lr E \lr 
  \Phi _{X\to X}^{\mathrsfs{S}}(E)\cong 0, $$
  we have $\dR \Hom (E, E)\otimes _{\mathbb{C}}E \cong E$. But again 
  this is impossible by the definition of $\mathbb{P}^n$-object. 
  $\quad \square$

\bibliographystyle{jplain}
\bibliography{math}

Yukinobu Toda, Graduate School of Mathematical Sciences, University of Tokyo

\textit{E-mail address}:toda@ms.u-tokyo.ac.jp

\end{document}